\documentclass[12]{article}
\usepackage[russian]{babel}
\usepackage{amsmath}
\usepackage {multirow}
\usepackage{amsfonts, amssymb, textcomp}
\input{tcilatex}
\begin{document}
\renewcommand{\baselinestretch}{1.6}
\newcounter{abc}
\textwidth = 161mm        
\textheight= 242mm        
\oddsidemargin= 11.7mm      
\evensidemargin= 1mm      
\topmargin=-5mm          
\setcounter{MaxMatrixCols}{10}
\mathsurround=2pt \sloppy \hfuzz=2pt \tolerance=6000
\emergencystretch=20pt \hyphenpenalty=10 \mathsurround=2pt
\righthyphenmin=2
\thispagestyle{empty}

\begin{center}
\textbf{Sabziyev N.M.}

\textbf{On a new method for calculation of the number of prime numbers in
the given interval} 
{nariman.sabziyev@gmail.com}
\end{center}

\textbf{Abstract.} Definition of the number of prime numbers in the interval.

\textbf{Key words.} Prime numbers.

For describing the behavior of the distribution function of prime numbers it
is necessary to investigate some auxiliary classes of natural numbers. The
given paper is devoted to investigation of some subclasses of natural
numbers, that allows to describe the distribution function of prime numbers
[1], [2].

\textbf{1. Denotation and necessary facts ([3], [4]).}

Denote by $N$ the set of all natural numbers, by $N\left( \leq n\right) $
the set of natural numbers not exceeding $n$.

It is obvious that the set $N\left( \geq 5\right) $ may be represented in
the form 
\[
N\left( \geq 5\right) =\left\{ 5,6,7,...,n,...\right\} =
\]%
\[
=\left\{ 6m-1,6m,6m+1,6m+2,6m+3,6m+4\right\} _{m=1}^{\infty }=
\]%
$$
=\left\{ 6m+\alpha ,\alpha =-1,0,1,2,3,4\right\} _{m=1}^{\infty }.\eqno(1.1)
$$

From (1.1) it is seen that the natural numbers of the form $6m$, $6m+2$, $%
6m+3$ and $6m+4$ for any natural value of $m$ are composite numbers; and
only natural numbers $6m-1$ and $6m+1$ $\left( m\in N\right) $ may be also
prime numbers.

Therefore we have [3]

\textbf{Theorem 1.1.} If $n$ is a prime natural number, then it is necessary
that $n$ should have the form $n=6m-1$ or $n=6m+1$.

If the number $6m+1$ (or $6m-1$) is a composite number, then at least it has
two cofactors, i.e.%
$$
6m\pm 1=\left( 6i+\alpha \right) \left( 6j+\beta \right) \eqno(1.2)
$$%
$$
\alpha \ \ \ \ and\ \ \ \ \beta =-1,0,1,2,3,4.\eqno(1.3)
$$

From (1.2) we have 
$$
6m\pm 1=6ij+\beta 6i+\alpha 6j+\alpha \beta =6m_{0}+\alpha \beta .\eqno(1.4)
$$

By the prime calculations from (1.4) we get 
$$
\alpha \beta =\pm 1,\ \ \ where\ \ \ \alpha \ \ \ and\ \ \ \beta \eqno(1.5)
$$%
were determined in (1.3).

Taking into account (1.5), from (1.2) we get 
$$
6m+1=\left( 6i\pm 1\right) \left( 6j\pm 1\right) \ \ \ \ i,j\in N,\ \ i\geq
j,\eqno\left( 1.6\right)
$$%
$$
6m-1=\left( 6i\pm 1\right) \left( 6j\mp 1\right) \ \ \ \ i,j\in N,\ \ i\geq
j.\eqno\left( 1.7\right)
$$

From (1.6) we get 
$$
m=6ij\pm \left( i+j\right) ,\ \ \ \left( i\geq j\right) \eqno\left(
1.8\right)
$$%
from (1.7) we have 
$$
m=6ij\pm \left( i-j\right) ,\ \ \ \left( i\geq j\right) .\eqno\left(
1.9\right)
$$

Denote the set of natural numbers of the form (1.6) and (1.7) by $M_{1}$ and 
$M_{2}$, respectively. Then we have

\textbf{Theorem 1.2.} For the natural number of the form $6m+1$ (or $6m-1$), 
$m\in N$ \ to be a composite natural number, it is necessary and sufficient
that $m\in M_{1}$ (or $m\in M_{2}$). 

Indeed, if $m\in M_{1}$ then $m=6ij\pm
\left( i+j\right) $. Hence $6m+1=36ij\pm 6\left( i+j\right) +1=\left( 6i\pm
1\right) \left( 6j\pm 1\right) $ and vice versa, if $6m+1$ is a composite
number, then $6m+1=\left( 6i\pm 1\right) \left( 6j\pm 1\right) $, hence it
follows that $m=6ij\pm \left( i+j\right) \in M$.

Denote by $H_{1}=N\backslash M_{1}\left( H_{2}=N\backslash M_{2}\right) $,
where $H_{1}\cap M_{1}\neq \varnothing $; $\left( H_{2}\cap M_{2}\right)
\neq \varnothing $.

Then we have

\textbf{Theorem 1.3.} For the natural number of the form $6m+1$ (or $6m-1$)
to be a prime natural number, it is necessary and sufficient that $m\in
H_{1} $ (or $m\in H_{2}$).

\textbf{2. The property of the set }$M_{1}\left( \leq m\right) $\textbf{.}

By definition, if $n\in M_{1}\left( \leq m\right) $, then $n=6ij\pm \left(
i+j\right) \leq m$, $i,j\in N$.

Note that if $m=\underset{i,j}{\max }\left( 6ij-i-j\right) $, then there
exists a natural number $\nu $ for $i=j=\nu $, $6\nu ^{2}-2\nu =m$, hence $%
\nu =\left[ \frac{1+\sqrt{6m+1}}{6}\right] $ i.e. $1\leq i\leq \nu _{j}$,
and if $m=\underset{i,j}{\max }\left( 6ij+i+j\right) $, then for $i=j=k$, $%
6k^{2}+2k=m$, hence $k=\left[ \frac{-1+\sqrt{6m+1}}{6}\right] $ and $1\leq
j\leq k$.

Denote by 
\[
A_{i}\left( m\right) =\left\{ 6ij-i-j\right\} =\left\{ \left( 6i-1\right)
j-i\leq m,\ \ \ i\leq j,\ \ i,j\in N\right\}
\]%
\[
B_{j}\left( m\right) =\left\{ 6ij+i+j\right\} =\left\{ \left( 6j+1\right)
j+i\leq m,\ \ \ i\leq j,\ \ i,j\in N\right\}
\]%
\[
1\leq i\leq \nu ,\ \ \ 1\leq j\leq k
\]
where $6i-1$ and $6j+1$ are only prime numbers.

Call $A_{i}\left( m\right) $ and $B_{i}\left( m\right) $ the subclasses with
prime coefficients of the set $M_{1}\left( \leq m\right) $.

It should be noted that 
\[
M_1\left( \leq m\right) =\left( \overset{V}{\underset{i=1}{\dbigcup }}%
A_{i}\left( m\right) \right) \cup \left( \overset{k}{\underset{i=1}{\dbigcup 
}}B_{j}\left( m\right) \right) .
\]

Denote the set of prime coefficients of the subclasses $A_{i}\left( m\right) $
and $B_{j}\left( m\right) $, by $K_{1}\left( -\right) $ and $K_{1}\left(
+\right) $:%
\[
K_{1}\left( -\right) =\left\{ 5,11,17,...,6\nu -1\right\} =\left\{
K_{1}^{\left( 1\right) }\left( -\right) ,K_{1}^{\left( 2\right) }\left(
-\right) ,...,K_{1}^{\left( \nu \right) }\left( -\right) \right\} ,
\]%
i.e. 
\[
K_{1}^{\left( 1\right) }\left( -\right) =5,\ \ K_{1}^{\left( 2\right)
}\left( -\right) =11,\ \ \ K_{1}^{\left( 3\right) }\left( -\right)
=17,...,K_{1}^{\left( \nu \right) }\left( -\right) =6\nu -1
\]
respectively, 
\[
K_{1}\left( +\right) =\left\{ 7,13,19,...,6k+1\right\} =\left\{
K_{1}^{\left( 1\right) }\left( +\right) ,K_{1}^{\left( 2\right) }\left(
+\right) ,...,K_{1}^{\left( {k} \right) }\left( +\right) \right\} ,
\]%
i.e.%
\[
K_{1}^{\left( 1\right) }\left( +\right) =7,\ \ K_{1}^{\left( 2\right)
}\left( +\right) =13,,...,\ \ \ K_{1}^{\left( k\right) }\left( -\right) =6k+1
\]%
where $6i-1$ and $\ 6j+1$ are only prime numbers and the number of the elements of the set $K_{2}\left( -\right) $ equals 
$$\nu _{2}\left( -\right) =C_{\nu _{0}}^{1}\cdot
C_{k_{0}}^{1}.$$

Denote by $K_{2}\left( -\right) $ the set with the elements of the form $%
6\tau -1\left( \tau \in N\right) $ being the product of two elements of the
set $K_{1}\left( -\right) \cup K_{1}\left( +\right) $:%
\[
K_{2}\left( -\right) =\left\{ 5\cdot 7,5\cdot 13,...,11\cdot 7,11\cdot
13,...,\left( 6\nu -1\right) \left( 6k+1\right) \right\} =
\]%
\[
=\left\{ K_{2}^{\left( 1\right) }\left( -\right) ,K_{2}^{\left( 2\right)
}\left( -\right) ,...,K_{2}^{\left( \gamma _{2}\left( -\right) \right)
}\left( -\right) \right\} ,
\]%
where $$K_{2}^{\left( 1\right) }\left( -\right) =5\cdot
7,K_{2}^{\left( 2\right) }\left( -\right) =5\cdot 13,...
$$
Here $\nu _{0}$ is the number of prime elements of the set $K_{1}\left( -\right) $, $k_{0}$ is the number of prime elements of the set $K_{1}\left( +\right) $. Denote by $ K_{2}\left( +\right) $ the set with the elements of the form $6\tau +1\left( \tau \in N\right) $ being the product of two elements of the set $K_{1}\left( -\right) \cup
K_{1}\left( +\right) $:%
\[
K_{2}\left( +\right) =\left\{ 5\cdot 11,5\cdot 17,...,7\cdot 13,7\cdot
19,...\right\} =
\]%
\[
=\left\{ K_{2}^{\left( 1\right) }\left( +\right) ,K_{2}^{\left( 2\right)
}\left( +\right) .,...,K_{2}^{\left( \gamma _{2}\left( +\right) \right)
}\right\} ,
\]%
where the number of the elements of the set $K_{2}\left( +\right) $ equals $%
\gamma _{2}\left( +\right) =V_{\nu _{0}}^{2}+C_{k_{0}}^{2}$%
\[
K_{2}\left( +\right) =5\cdot 11,K_{2}^{\left( 2\right) }\left( +\right)
=5\cdot 17,...,K_{2}^{\left( 3\right) }\left( +\right) =7\cdot 13,...
\]

The set $\ K_{q}\left( -\right) $ and $K_{q}\left( +\right) $, where $2\leq
i+j=q$ , \ $1\leq i\leq \nu $ , $1\leq j\leq k$\ \ is determined in the same
way.

Now we can calculate the number of the elements of subclasses of the set $%
M_{1}\left( \leq m\right) $:%
\[
mesA_{1}\left( \leq m\right) =mes\left\{ 5t-1\leq m\right\} =\left[ \frac{m+1%
}{5}\right] =
\]%
$$
=\left[ \frac{m+\frac{K_{1}^{\left( 1\right) }\left( -\right) +1}{6}}{%
K_{1}^{\left( 1\right) }\left( -\right) }\right] =\left[ \frac{%
6m+K_{1}^{\left( 1\right) }\left( -\right) +1}{6K_{1}^{1}\left( -\right) }%
\right] \eqno\left( 2.1\right)
$$%
\[
mesA_{2}\left( \leq m\right) =mes\left\{ 11t-2\leq m\right\} =\left[ \frac{%
m+2}{11}\right] =
\]%
$$
=\left[ \frac{m+\frac{K_{1}^{\left( 2\right) }\left( -\right) +1}{6}}{%
K_{1}^{\left( 1\right) }\left( -\right) }\right] =\left[ \frac{%
6m+K_{1}^{\left( 2\right) }\left( -\right) +1}{6K_{1}^{2}\left( -\right) }%
\right] \eqno\left( 2.2\right)$$
\begin {center}
............................................................\\
............................................................
\end {center}
\[mesA_{\nu }\left( \leq m\right) =mes\left\{ \left( 6\nu -1\right) t-\nu \leq
m\right\} =\left[ \frac{m+\nu }{6\nu -1}\right] =
\]%
$$
=\left[ \frac{6m+K_{1}^{\left( \nu \right) }\left( -\right) +1}{6K_{1}^{\nu
}\left( -\right) }\right] \eqno\left( 2.3\right)
$$%
where $K_{1}^{\left( 1\right) }\left( -\right) $, $K_{1}^{\left( 2\right)
}\left( -\right) $, ... , $K_{1}^{\left( \nu \right) }\left( -\right) $ are
the elements of the set $K_{1}\left( -\right) $; 
\[
mesB_{1}\left( \leq m\right) =mes\left\{ 7t+1\leq m\right\} =\left[ \frac{m-1%
}{7}\right] =
\]%
$$
=\left[ \frac{6m-K_{1}^{\left( 1\right) }\left( -\right) +1}{6K_{1}^{\nu
}\left( +\right) }\right] \eqno\left( 2.4\right)
$$%
\[
mesB_{2}\left( \leq m\right) =mes\left\{ 13t+2\leq m\right\} =\left[ \frac{%
m-2}{13}\right] =
\]%
$$
=\left[ \frac{6m-K_{1}^{\left( 2\right) }\left( +\right) +1}{%
6K_{1}^{2}\left( +\right) }\right] \eqno\left( 2.5\right)
$$%
\begin {center}
............................................................\\
............................................................
\end {center}
\[
mesB_{k}\left( \leq m\right) =mes\left\{ \left( 6k+1\right) t+k\right\} = 
\left[ \frac{m-k}{6k+1}\right] =
\]%
$$
=\left[ \frac{6m-K_{1}^{\left( k\right) }\left( +\right) +1}{6K_{1}^{\left(
k\right) }\left( +\right) }\right] ,\eqno\left( 2.6\right)
$$%
where $K_{1}^{\left( 1\right) }\left( +\right) $, $K_{1}^{\left( 2\right)
}\left( +\right) $,...,$K_{1}^{\left( k\right) }\left( +\right) $ are the
elements of the set $K_{1}\left( +\right) $.

Let $m\in A_{1}\left( m\right) \cap A_{2}\left( m\right) $. Then $%
m=5t-1=11\tau -2$, $5t=11\tau -1$, $t=2\tau +\frac{\tau -1}{5}$, $\tau
=5\tau _{1}+1$, $m=5\cdot 11\tau _{1}+9$ and $m=5\cdot 11\tau _{0}-46\left(
\tau =\tau _{1}-1\right) $. Hence it is seen that the coefficient $\tau l0$
represents the number of the form $6t+1$, i.e. $5\cdot 11$ is a natural
number of the form $6t+1$.

The number of the elements of the set $A_{1}\left( m\right) \cap A_{2}\left(
m\right) $ equals 
$$
\left[ \frac{m+46}{5\cdot 11}\right] =\left[ \frac{m+\frac{5\cdot 5\cdot 11+1%
}{6}}{5\cdot 11}\right] =\left[ \frac{6m+5\cdot 5\cdot 11+1}{6\cdot 5\cdot 11%
}\right] .\eqno\left( 2.7\right)
$$

And if $m\in A_{1}\left( m\right) \cap B_{1}\left( m\right) $, then 
\[
5t-1=7\tau +1,\ \ \ 5t=7\tau +2,\ \ \ t=\tau +2\frac{\tau +1}{5},\ \ \ \tau
=5\tau _{1}-1,\ \ \ m=5\cdot 7\tau _{1}-6,
\]
hence it is seen that the coefficient $\tau _{1}$ represents the number of
the form $6t-1$, i.e. 5.7 is a natural number of the form $6t-1$; the number
of the elements of the set $A_{1}\left( m\right) \cap B_{1}\left( m\right) $
equals 
\[
\left[ \frac{m+6}{5\cdot 7}\right] =\left[ \frac{m+\frac{5\cdot 7+1}{6}}{%
5\cdot 7}\right] =\left[ \frac{6m+5\cdot 7+1}{6\cdot 5\cdot 7}\right] .
\]

Similarly continuing (by induction), we find that the number of the elements
of the set 
\[
s\left( m\right) =\left( \underset{s_{1}=0}{\overset{s}{\dbigcap }}%
A_{i_{s_{1}}}\left( m\right) \right) \cap \left( \underset{r_{1}=0}{\overset{%
r}{\dbigcap }}B_{j_{r_{1}}}\left( m\right) \right) ,\ \ \ 2\leq
s_{1}+r_{1}\leq q,\ \ 2\leq q\leq \nu +k
\]%
we get 
$$
mes S\left( m\right) =\left[ \frac{6m+a\overset{s}{\underset{s_{1}=0}{\dprod }%
}\left( 6s_{1}-1\right) \overset{r}{\underset{r_{1}=0}{\dprod }}\left(
6r_{1}+1\right) +1}{6\overset{s}{\underset{s_{1}=0}{\dprod }}\left(
6s_{1}-1\right) \overset{r}{\underset{r_{1}=0}{\dprod }}\left(
6r_{1}+1\right) }\right] ,\eqno(2.8)
$$%
where 
\[
a=\left\{ 
\begin{array}{c}
1,\ \ \ if\ \ s\ \ is\ an\ add\ number; \\ 
5\ \ if\ \ s\ \ is\ an\ even\ number,%
\end{array}%
\right.
\]%
\[
a\overset{0}{\underset{s_{1}=0}{\dprod }}\left( -1\right) =1,\ \ \ 1\leq
i_{s_{1}}\leq \nu ,\ \ \ 1\leq j_{r_{1}}\leq k.
\]

If we denote the number of the composite numbers of the form $6t+1\left(
t\in N\right) $\ by $P^{\left( +\right) }\left( 6m+1\right) $, then from
(2.1)-(2.8) we have

\textbf{Theorem 2.1.} For the given $m\in N$ the number of the elements of
the set $M_{1}\left( \leq m\right) $ (i.e. the number of the composite
numbers of the form $6\tau +1\left( \tau \in N\right) $ doesn't exceed $6m+1$%
) equals 
\[
P^{\left( +\right) }\left( 6m+1\right) =\underset{i=1}{\overset{\nu }{\sum }}%
\left[ \frac{6m+K_{1}^{\left( i\right) }\left( -\right) +1}{6K_{1}^{\left(
i\right) }\left( -\right) }\right] +\underset{j=1}{\overset{k}{\sum }}\left[ 
\frac{6m-K_{1}^{\left( j\right) }\left( +\right) +1}{6K_{1}^{\left( j\right)
}\left( +\right) }\right] +
\]%
\[
+\underset{q=2}{\overset{\nu +k}{\sum }}\left( -1\right) ^{q-1}\left( 
\underset{i=1}{\overset{\gamma _{q}\left( -\right) }{\sum }}\left[ \frac{%
6m+K_{q}^{\left( i\right) }\left( -\right) +1}{6K_{q}^{\left( i\right)
}\left( -\right) }\right] +\underset{j=1}{\overset{\gamma _{q}\left(
+\right) }{\sum }}\left[ \frac{6m+5K_{q}^{\left( j\right) }\left( +\right) +1%
}{6K_{q}^{\left( j\right) }\left( +\right) }\right] \right)
\]%
where 
\[
\gamma _{2}\left( -\right) =C_{\nu }^{1},C_{k}^{1},\gamma _{3}\left(
-\right) =C_{\nu }^{1}C_{k}^{2}+C_{\nu }^{3},...
\]%
\[
\gamma _{2}\left( +\right) =C_{\nu }^{2}+C_{k}^{2},\ \ \gamma _{3}\left(
+\right) =C_{\nu }^{2}C_{k}^{1}+C_{k}^{3},...
\]

Denote by $\pi ^{\left( +\right) }\left( 6m+1\right) $ the number of prime
numbers of the form $6\tau +1\left( \tau \in N\right) $, then by $%
H_{1}\left( \leq m\right) =N\left( \leq m\right) \backslash M_{1}\left( \leq
m\right) $ it holds

\textbf{Theorem 2.2.} For the given $m$, the number of prime numbers of the
form $6\tau +1\left( \tau \in N\right) $ not exceeding $6m+1$ (i.e. the
number of the elements of the set $H_{1}\left( \leq m\right) $) equals 
$$
\pi ^{\left( +\right) }\left( 6m+1\right) =m-P^{\left( +\right) }\eqno( 2.10)
$$%
where $P^{\left( +\right) }\left( 6m+1\right) $ was determined in equality
(2.9), 
\[
\nu =\left[ \frac{1+\sqrt{6m+1}}{6}\right] ,\ \ \ k=\left[ \frac{-1+\sqrt{%
6m+1}}{6}\right] .
\]

\textbf{Example 1.} Let $m=50$, then $\nu =3$, $k=2$%
\[
K_{1}\left( -\right) =\left\{ 5,11,17\right\} =\left\{ K_{1}^{\left(
1\right) }\left( -\right) ,\ K_{1}^{\left( 2\right) }\left( -\right) ,\
K_{1}^{\left( 3\right) }\left( -\right) \right\} ,\ \ i.e.
\]%
\[
K_{1}^{\left( 1\right) }\left( -\right) =5,\ \ K_{1}^{\left( 2\right)
}\left( -\right) =11,\ \ \ K_{1}^{\left( 3\right) }\left( -\right) =17.
\]%
\[
K_{1}\left( +\right) =\left\{ 7,13\right\} =\left\{ K_{1}^{\left( 1\right)
}\left( +\right) ,\ K_{1}^{\left( 2\right) }\left( +\right) \right\} ,\ \
i.e.
\]%
\[
K_{1}^{\left( 1\right) }\left( +\right) =7,\ \ K_{1}^{\left( 2\right)
}\left( -\right) =13.
\]

Then 
\[
\overset{3}{\underset{i=1}{\sum }}\left[ \frac{m+\frac{K_{1}^{\left(
i\right) }\left( -\right) +1}{6}}{K_{1}^{\left( i\right) }\left( -\right) }%
\right] +\overset{2}{\underset{j=1}{\sum }}\left[ \frac{m-\frac{%
K_{1}^{\left( j\right) }\left( +\right) -1}{6}}{K_{1}^{\left( j\right)
}\left( +\right) }\right] =
\]%
$$
=\left[ \frac{m+1}{5}\right] +\left[ \frac{m+2}{11}\right] +\left[ \frac{m+3%
}{17}\right] +\left[ \frac{m-1}{7}\right] +\left[ \frac{m-2}{13}\right] .%
\eqno\left( 2.11\right)\\$$%

Since
\[
K_{2}\left( -\right) =\left\{ 5\cdot 7,5\cdot 13,11\cdot 7,11\cdot
13,17\cdot 7,17\cdot 13\right\} .
\]

The elements $K_{2}\left( -\right) $ have the form $6\tau -1$, and the
number of the elements of the set $K_{2}\left( -\right) $ equals $\gamma
_{2}\left( -\right) =C_{3}^{1}\cdot C_{2}^{1}=6$, then 
\[
\overset{\gamma _{2}\left( -\right) }{\underset{i=1}{\sum }}\left[ \frac{m+%
\frac{K_{2}^{\left( i\right) }\left( -\right) +1}{6}}{K_{2}^{\left( i\right)
}\left( -\right) }\right] =\left[ \frac{m+6}{5\cdot 7}\right] +\left[ \frac{%
m+11}{5\cdot 13}\right] +
\]%
$$
+\left[ \frac{m+13}{11\cdot 7}\right] +\left[ \frac{m+24}{11\cdot 13}\right]
+\left[ \frac{m+20}{7\cdot 11}\right] +\left[ \frac{m37}{13\cdot 17}\right] .%
\eqno\left( 2.12\right)
$$%
and%
\[
K_{2}\left( +\right) =\left\{ 5\cdot 11,5\cdot 17,11\cdot 17,7\cdot
13\right\} ,
\]%
the elements $K_{2}\left( +\right) $ have the form $6\tau +1$, and the
number of the elements of the set $K_{2}\left( +\right) $ equals $\gamma
_{2}\left( +\right) =C_{\nu }^{2}+C_{k}^{2}=4$ 
\[
\overset{\gamma _{2}\left( +\right) }{\underset{i=1}{\sum }}\left[ \frac{m+%
\frac{K_{2}^{\left( i\right) }\left( +\right) +1}{6}}{K_{2}^{\left( i\right)
}\left( +\right) }\right] =\left[ \frac{m+46}{5\cdot 11}\right] +\left[ 
\frac{m+71}{5\cdot 17}\right] +
\]%
$$
+\left[ \frac{m+156}{11\cdot 17}\right] +\left[ \frac{m+76}{7\cdot 13}\right]
.\eqno\left( 2.13\right)
$$%
Further,
\[
K_{3}\left( -\right) =\left\{ 5\cdot 11\cdot 17,5\cdot 7\cdot 13,11\cdot
7\cdot 13,17\cdot 7\cdot 13\right\} ,
\]%
whose elements have the form $6\tau -1$, and the number of the elements $%
K_{3}\left( -\right) $ equals $\gamma _{3}\left( -\right) =C_{\nu
}^{3}+C_{\nu }^{1}\cdot C_{k}^{2}=4$%
\[
\overset{\gamma _{3}\left( -\right) }{\underset{i=1}{\sum }}\left[ \frac{m+%
\frac{K_{3}^{\left( i\right) }\left( -\right) +1}{6}}{K_{3}^{\left( i\right)
}\left( -\right) }\right] =\left[ \frac{m+156}{5\cdot 11\cdot 17}\right] +%
\left[ \frac{m+76}{5\cdot 7\cdot 13}\right] +
\]%
$$
+\left[ \frac{m+167}{7\cdot 11\cdot 13}\right] +\left[ \frac{m+258}{7\cdot
13\cdot 17}\right] .\eqno\left( 2.14\right)
$$%

Since
\[
K_{3}\left( +\right) =\left\{ 5\cdot 11\cdot 7,5\cdot 17\cdot 7,11\cdot
17\cdot 7,5\cdot 11\cdot 13,5\cdot 17\cdot 13,11\cdot 17\cdot 13\right\} .
\]

The elements $K_{3}\left( +\right) $ have the form $6\tau +1$, and the
number of the elements equals $\gamma _{3}\left( +\right) =C_{\nu }^{2}\cdot
C_{k}^{1}=6$%
\[
\overset{\gamma _{3}\left( +\right) }{\underset{j=1}{\sum }}\left[ \frac{m+%
\frac{5K_{3}^{\left( j\right) }\left( +\right) +1}{6}}{K_{3}^{\left(
j\right) }\left( +\right) }\right] =\left[ \frac{m+321}{5\cdot 7\cdot 11}%
\right] +\left[ \frac{m+496}{5\cdot 7\cdot 17}\right] +
\]%
$$
+\left[ \frac{m+1091}{7\cdot 11\cdot 13}\right] +\left[ \frac{m+596}{5\cdot
11\cdot 13}\right] +\left[ \frac{m+921}{5\cdot 13\cdot 17}\right] +\left[ 
\frac{m+2026}{11\cdot 13\cdot 17}\right] .\eqno\left( 2.15\right)
$$%
For
\[
K_{4}\left( -\right) =\left\{ 5\cdot 11\cdot 17\cdot 7,\ \ 5\cdot 11\cdot
17\cdot 13\right\} ,
\]%
the elements $K_{4}\left( -\right) $ have the form $6\tau -1$, and the
number of the elements equals $\gamma _{4}\left( -\right) =C_{\nu }^{\nu
}\cdot C_{k}^{1}=2$%
$$
\overset{\gamma _{4}\left( -\right) }{\underset{i=1}{\sum }}\left[ \frac{m+%
\frac{K_{4}^{\left( j\right) }\left( -\right) +1}{6}}{K_{4}^{\left( i\right)
}\left( -\right) }\right] =\left[ \frac{m+1091}{5\cdot 7\cdot 11\cdot 13}%
\right] +\left[ \frac{m+2026}{5\cdot 11\cdot 13\cdot 17}\right] \eqno\left(
2.16\right)$$%
for
\[
K_{4}\left( +\right) =\left\{ 5\cdot 11\cdot 7\cdot 13,\ \ 5\cdot 17\cdot
7\cdot 13,\ 11\cdot 17\cdot 7\cdot 13\right\} ,
\]%
the elements $K_{4}\left( +\right) $ are of the form $6\tau +1$, and the
number of the elements equals $\gamma _{4}\left( +\right) =C_{\nu }^{2}\cdot
C_{k}^{k}=3$%
\[
\overset{\gamma _{4}\left( +\right) }{\underset{j=1}{\sum }}\left[ \frac{m+%
\frac{5K_{4}^{\left( j\right) }\left( +\right) +1}{6}}{K_{4}^{\left(
j\right) }\left( +\right) }\right] =\left[ \frac{m+4171}{5\cdot 7\cdot
11\cdot 13}\right] +
\]%
$$
+\left[ \frac{m+6446}{5\cdot 7\cdot 13\cdot 17}\right] +\left[ \frac{m+1418}{%
7\cdot 11\cdot 13\cdot 17}\right] .\eqno\left( 2.17\right)
$$%
\[
K_{5}\left( -\right) =\left\{ 5\cdot 11\cdot 17\cdot 7\cdot 13\right\} ,\ \
\gamma _{5}\left( -\right) =C_{\nu }^{\nu }\cdot C_{k}^{k}=1
\]%
and the number of the elements equals 
$$
\left[ \frac{m+14181}{5\cdot 7\cdot 11\cdot 13\cdot 17}\right] \eqno\left(
2.18\right)
$$%
and 
\[
K_{5}\left( +\right) \equiv \varnothing .
\]

Thus, taking into account equalities (2.11)-(2.18), we get 
$$
P^{\left( +\right) }\left( 301\right) =22\eqno\left( 2.19\right)
$$%
and from (2.10) we have 
$$
\pi ^{\left( +\right) }\left( 6\cdot 50+1\right) =\pi ^{\left( +\right)
}\left( 306\right) =50-22=28.\eqno\left( 2.20\right)
$$

\textbf{3. Property of the set }$M_{2}\left( \leq m\right) $.

By definition, if $n\in M_{2}\left( \leq m\right) $ then 
\[
n=6it-i+t\leq m,\ \ \ i\leq t,\ \ i,t\in N
\]%
or%
\[
n=6jt+j-t\leq m,\ \ \ j\leq t,\ \ j,t\in N.
\]%
where $m=\underset{i,t}{\max }\left\{ 6it-i+t\right\} $ or $m=\underset{j,t}{%
\max }\left\{ 6jt+j-t\right\} $. Then there exists a natural number $r$, for 
$i=j=t=r$ we have $6r^{2}=m$, hence we get 
\[
r=\left[ \frac{\sqrt{6m}}{6}\right] ,
\]
where $1\leq i\leq r$, $1\leq j\leq r$.

Denote by 
\[
C_{i}\left( m\right) =\left\{ \left( 6i-1\right) t+i\leq m,\ \ \ i\leq t,\ \
i,t\in N,\ \ 1\leq i\leq r\right\}
\]
and 
\[
D_{j}\left( m\right) =\left\{ \left( 6j+1\right) t-i\leq m,\ \ \ j\leq t,\ \
j,t\in N,\ \ 1\leq j\leq r\right\}
\]
where $6i-1$ and $6j+1$ are only prime numbers.

Call $C_{i}\left( m\right) $ and $D_{j}\left( m\right) $ the subclasses with
prime coefficients of the set $M_{2}\left( \leq m\right) $. Obviously, 
\[
M_{2}\left( \leq m\right) =\left( \overset{r}{\underset{i=1}{\dbigcup }}%
C_{i}\left( m\right) \right) \cup \left( \overset{r}{\underset{j=1}{\dbigcup 
}}D_{j}\left( m\right) \right) .
\]

Denote by $K_{1}^{\left( 1\right) }\left( -\right) $ and $K_{2}^{\left(
1\right) }\left( +\right) $ the set of prime coefficients of the subclasses $%
C_{i}\left( m\right) $ and $D_{j}\left( m\right) $\ and the set $M_{2}\left(
\leq m\right) $:%
\[
K_{1}\left( -\right) =\left\{ 5,11,17,...,6r-1\right\} ,
\]%
\[
K_{1}\left( +\right) =\left\{ 7,13,19,...,6r+1\right\} ,
\]
where the elements of the
set $K_{1}\left( -\right) \cup K_{1}\left( +\right) $ are only prime numbers.

As in the set $M_{1}\left( \leq m\right) $, here we also determine the set 
\[
K_{2}\left( -\right) ,\ \ K_{2}\left( +\right) ,\ \ K_{3}\left( -\right) ,\
\ K_{3}\left( +\right) ,...,K_{q}\left( -\right) \ \ and\ \ \ K_{q}\left(
+\right) ,
\]
and calculate the number of the elements of the subclasses $C_{i}\left(
m\right) $ and $D_{j}\left( m\right) $%
\[
mes\left( C_{1}\left( \leq m\right) \right) =mes\left( 5t+1\leq m\right) = 
\left[ \frac{m-1}{5}\right] =\left[ \frac{6m-K_{1}^{\left( 1\right) }\left(
-\right) -1}{6K_{1}^{\left( 1\right) }\left( -\right) }\right] ,
\]%
\[
mes\left( C_{2}\left( \leq m\right) \right) =mes\left( 11t+2\leq m\right) = 
\left[ \frac{m-2}{11}\right] =\left[ \frac{6m-K_{1}^{\left( 2\right) }\left(
-\right) -1}{6K_{1}^{\left( 2\right) }\left( -\right) }\right] ,
\]%
\[
.................................................................................
\]%
\[
.................................................................................
\]%
\[
mes\left( D_{1}\left( \leq m\right) \right) =mes\left( 7t-1\leq m\right) = 
\left[ \frac{m+1}{7}\right] =\left[ \frac{6m+K_{1}^{\left( 1\right) }\left(
+\right) -1}{6m}\right] ,
\]%
\[
mes\left( D_{2}\left( \leq m\right) \right) =mes\left( 13t-2\leq m\right) = 
\left[ \frac{m+2}{13}\right] =\left[ \frac{6m+K_{1}^{\left( 2\right) }\left(
+\right) -1}{6m}\right] ,
\]
and etc.

If 
\[
R\left( m\right) =\left( \underset{s_{1}=0}{\overset{s}{\dbigcap }}%
C_{i_{S_{1}}}\left( m\right) \right) \cap \left( \underset{r_{1}=0}{\overset{%
r}{\dbigcap }}D_{jr_{1}}\right) ,
\]%
\[
2\leq s_{1}+r_{1}=q\leq s+r
\]
then 
\[
mesR\left( m\right) =\left[ \frac{6m+b\overset{s}{\underset{s_{1}=0}{\dprod }%
}\left( 6s_{1}-1\right) \overset{r}{\underset{r_{1}=0}{\dprod }}\left(
6r_{1}+1\right) }{6\overset{s}{\underset{s_{1}=0}{\dprod }}\left(
6s_{1}-1\right) \overset{r}{\underset{r_{1}=0}{\dprod }}\left(
6r_{1}+1\right) }\right]
\]
where 
\[
b=\left\{ 
\begin{array}{c}
1,\ \ \ if\ \ \ \ s\ \ is\ an\ even\ number \\ 
5,\ \ \ if\ \ s\ \ is\ an\ odd\ number%
\end{array}%
\right.
\]%
\[
\overset{0}{\underset{0}{\dprod }}\left( -1\right) =1,\ \ \ 1\leq i,\ \ \
j\leq r
\]%
where $K_{q}^{\left( 1\right) }\left( +\right) $, $K_{q}^{\left( 1\right)
}\left( -\right) $, $K_{q}^{\left( 2\right) }\left( +\right) $, $%
K_{q}^{\left( 2\right) }\left( -\right) $,... is determined as in
calculating the number of the elements of subclasses of the set $M_{1}\left(
\leq m\right) $.

Denote by $P^{\left( -\right) }\left( 6m-1\right) $ the number of composite
numbers of the form $6\tau -1$ $\left( \tau \in N\right) $ not exceeding $%
6m-1$, then we have

\textbf{Theorem 3.1.} For the given $m\in N$, the number of the elements of
the set $M_{2}\left( \leq m\right) $ (i.e. the number of composite numbers
of the form $6\tau -1\left( \tau \in N\right) $ not exceeding $6m-1$) equals 
\[
P^{\left( -\right) }\left( 6m-1\right) =\overset{r}{\underset{i=1}{\sum }}%
\left[ \frac{6m+K_{1}^{\left( i\right) }\left( -\right) +1}{6K_{1}^{\left(
i\right) }\left( -\right) }\right] +\overset{r}{\underset{j=1}{\sum }}\left[ 
\frac{6m-K_{1}^{\left( j\right) }\left( +\right) +1}{6K_{1}^{\left( j\right)
}\left( +\right) }\right] +
\]%
\[
+\overset{r+r}{\underset{q=2}{\sum }}\left( -1\right) ^{q-1}\left( \overset{%
\gamma _{q}\left( -\right) }{\underset{i=1}{\sum }}\left[ \frac{%
6m+5K_{q}^{\left( i\right) }\left( -\right) +1}{6K_{q}^{\left( i\right)
}\left( -\right) }\right] +\overset{\gamma _{q}\left( +\right) }{\underset{%
j=1}{\sum }}\left[ \frac{6m+K_{q}^{\left( j\right) }\left( +\right) -1}{%
6K_{q}^{\left( j\right) }\left( +\right) }\right] \right)
\]
where $\gamma _{q}\left( -\right) $ and $\gamma _{q}\left( +\right) $\ is
determined as in theorem 2.1.

Denote by $\pi ^{\left( -\right) }\left( 6m-1\right) $ the number of prime
numbers of the form $6\tau -1$ $\left( \tau \in N\right) $ not exceeding $%
6m-1$, then by 
\[
H_{2}\left( \leq m\right) =N\left( \leq m\right) \backslash M_{2}\left( \leq
m\right)
\]
it holds

\textbf{Theorem 3.2.} For the given $m\in N$ the number of prime numbers of
the form $6\tau -1$ $\left( \tau \in N\right) $ not exceeding $6m-1$ (i.e.
the number of the elements of the set $H_{2}\left( \leq m\right) $) equals 
$$
\pi ^{\left( -\right) }\left( 6m-1\right) =m-P^{\left( -\right) }\left(
6m-1\right), \eqno\left( 3.2\right)
$$%
where $P^{\left( -\right) }\left( 6m-1\right) $ is determined as equality
(3.1), $r=\left[ \frac{\sqrt{6m}}{6}\right] $.

\textbf{Example 2.} Let $m=50$, then $r=2$ and 
\[
K\left( -\right) =\left\{ 5,11\right\} ,\ \ K\left( +\right) =\left\{
7,13\right\}
\]
i.e. 
\[
K_{1}^{1}\left( -\right) =5,\ \ K_{1}^{\left( 2\right) }\left( -\right)=11 ,\ \
K_{1}^{\left( 1\right) }\left( +\right) =7,\ \ K_{1}^{\left( 2\right)
}\left( +\right) =13.
\]
Then from (3.1) we have 
\[
P^{\left( -\right) }\left( 301\right) =\left[ \frac{50-1}{5}\right] +\left[ 
\frac{50-2}{11}\right] +\left[ \frac{50+1}{7}\right] +\left[ \frac{50+2}{13}%
\right] -
\]%
\[
\left( -\left[ \frac{50+29}{5\cdot 7}\right] +\left[ \frac{50+54}{5\cdot 13}%
\right] +\left[ \frac{50+64}{7\cdot 11}\right] +\left[ \frac{50+119}{11\cdot
13}\right] +\right.
\]%
\[
\left. +\left[ \frac{50+9}{5\cdot 11}\right] +\left[ \frac{50+15}{7\cdot 13}%
\right] \right) +\left( \left[ \frac{50+327}{5\cdot 7\cdot 13}\right] +\left[
\frac{50+834}{11\cdot 7\cdot 13}\right] +\right.
\]%
\[
\left. +\left[ \frac{50+64}{5\cdot 11\cdot 7}\right] +\left[ \frac{50+119}{%
5\cdot 11\cdot 13}\right] \right) -\left[ \frac{50+834}{5\cdot 7\cdot
11\cdot 13}\right] =
\]%
\[
=\left( 9+4+7+4\right) -\left( 2+1+1+1+1\right) =24-6=18
\]
\[
P^{\left( -1\right) }\left( 6m-1\right) =P^{\left( -\right) }\left(
301\right) =18
\]%
i.e.%
\[
\pi ^{\left( -\right) }\left( 6m-1\right) =\pi ^{\left( -\right) }\left(
299\right) =50-18=32.
\]

\textbf{4. Calculation of the number of prime numbers not exceeding }$6m+1$%
\textbf{.}

Denote by $\pi \left( 6m+1\right) $ the number of prime numbers not
exceeding $6m+1$. Then from theorems 2.2 and 3.2 we have

\textbf{Theorem 4.} The number of prime numbers not exceeding $6m+1$ (except
2 and 3) equals 
\[
\pi \left( 6m+1\right) =\pi ^{\left( +\right) }\left( 6m+1\right) +\pi
^{\left( -\right) }\left( 6m+1\right) =
\]%
\[
=2m-\left( P^{\left( +\right) }\left( 6m+1\right) +P^{\left( -\right)
}\left( 6m-1\right) \right) .
\]

\begin{center}
\bigskip

\bigskip

\textbf{References}
\end{center}

1. Bookhstab A.A. Theory of numbers. M. Nauka, 1966.

2. Ingam A.I. The distribution of prime numbers. M.L. 1936, 159 p.

3. Sabziyev N.M. Characteristics of sets wich are generating numbers of a
kind $6m\pm 1$. Transaction of the National Academy of Sciences of
Azerbaijan, 2003, No 3, pp. 41-49.

4. Sabziyev N.M. Distribution of prime numbers in natural rows. Transactions
of the National Azademy of Sciences of Azerbaijan, 2003, No 3, pp. 50-56.

\bigskip 

\end{document}